# Some factors of nonsingular Bernoulli shifts

by

Zemer Kosloff (Jerusalem) and Terry Soo (London)

**Abstract.** We give elementary constructions of factors of nonsingular Bernoulli shifts. In particular, we show that all nonsingular Bernoulli shifts on a finite set of symbols which satisfy the Doeblin condition have a factor that is equivalent to an independent and identically distributed system. We also prove that there are type-III$_1$ Bernoulli shifts of every possible ergodic index, answering a question of Danilenko and Lemańczyk [Ergodic Theory Dynam. Systems 39 (2019), 3292–3321].

**1. Introduction.** Let $A$ be a subset of $\mathbb{R}$, which will usually be an interval or a finite set. Let $(\rho_i)_{i \in \mathbb{Z}}$ be a sequence of probability measures on $A$. Let $\rho = \bigotimes_{i \in \mathbb{Z}} \rho_i$ be the product measure on $A^{\mathbb{Z}}$ and $T : A^{\mathbb{Z}} \to A^{\mathbb{Z}}$ be the *left-shift* given by $(Tx)_i = x_{i+1}$. We say that the dynamical system $(A^{\mathbb{Z}}, \mathcal{B}, \rho, T)$ is a *Bernoulli shift*, where $\mathcal{B}$ is the usual Borel product sigma-algebra. We say that the product measure $\rho$ and the corresponding system are *nonsingular* if the measure $\rho \circ T^{-1}$ is equivalent to $\rho$; such systems can be thought of as models for systems that are *not* in equilibrium. We say that the system is *conservative* if for all $E \in \mathcal{B}$ with $\rho(E) > 0$, there exists a nonzero integer $n$ such that $\rho(E \cap T^{-n}E) > 0$. If the measures $\rho_i$ are all identical, then the Bernoulli shift is nonsingular and conservative, and we say it is an *independent and identically distributed* (i.i.d.) system. Ornstein [32] proved that entropy is a complete isomorphism invariant for i.i.d. systems, but the case of nonsingular systems appears to be more delicate and it is unclear what role entropy plays [8, Section 9].

Let $\rho$ and $\nu$ be nonsingular product measures on $A^{\mathbb{Z}}$. We say that a measurable map $\phi : A^{\mathbb{Z}} \to A^{\mathbb{Z}}$ is a *factor* from $\rho$ to $\nu$ if $\phi$ is *equivariant* so that $\phi \circ T = T \circ \phi$ and the push-forward of $\rho$ under $\phi$ given by $\rho \circ \phi^{-1}$ is equivalent to $\nu$; in the case that the push-forward is $\nu$, we say that the factor is *measure-preserving*. If the inverse map $\phi^{-1}$ also serves as a factor from $\nu$ to $\rho$, then $\phi$ is an *isomorphism*, and we say that the corresponding Bernoulli shifts are *isomorphic*.

A product measure $\rho$ on $A^{\mathbb{Z}}$ satisfies the *Doeblin condition* if there exists $\delta > 0$ such that for all $a \in A$ and $i \in \mathbb{Z}$, we have $\rho_i(a) > \delta$.

THEOREM 1 (Low entropy i.i.d. factor). *Every nonsingular Bernoulli shift on a finite set of symbols that is equipped with a product measure which satisfies the Doeblin condition has an i.i.d. factor.*

Theorem 1 is a weak version of Sinai's celebrated factor theorem [37], for the nonsingular setting, that will apply to Bernoulli shifts considered by Krengel [26], Hamachi [16], Kosloff [21, 22], and Vaes and Wahl [40]; see formula (1).

An isomorphism invariant that is often considered in the study of nonsingular systems is orbit equivalence; two systems are *orbit equivalent* if there is a measurable bijection $\phi : A^{\mathbb{Z}} \to A^{\mathbb{Z}}$ such that the push-forward of $\rho$ is equivalent to $\nu$, and $\phi(\text{orb}(x)) = \text{orb}(\phi(x))$ for $\rho$-almost every $x \in A^{\mathbb{Z}}$, where $\text{orb}(x) = \{T^n(x) : n \in \mathbb{Z}\}$. Dye's theorem [10, 11] states that all ergodic nonatomic probability-preserving, and thus i.i.d., systems are orbit equivalent. A nonsingular and conservative Bernoulli shift with an associated product measure $\rho$ that *cannot* be endowed with a possible infinite shift-invariant measure that is equivalent to $\rho$ is said to be of *Krieger type*-III.

In this paper, we are particularly interested in constructing explicit factors between various type-III Bernoulli shifts and producing i.i.d. factors from these type-III shifts. We remark that although some of the results are confined to specific examples, it is nontrivial to prove the existence of a type-III Bernoulli shift, and the first construction of this type is due to Hamachi [16].

Krieger's theory [27, 28] further assigns a parameter $\lambda \in [0,1]$ to each type-III system and tells us that any two type-III$_\lambda$ systems are orbit equivalent, provided that $\lambda > 0$. Recently, in [25] we constructed the first type-III$_\lambda$ Bernoulli shifts for $\lambda \in (0,1)$, which were given by an explicit sequence of step functions on an interval taking the three values $1, \lambda$, and $\lambda^{-1}$. We will specify this construction in Section 3.2.

THEOREM 2. *Let $\lambda, \lambda' \in (0,1]$. There exists a type-III$_\lambda$ Bernoulli shift, specified in Section 3, which has a type-III$_{\lambda'}$ Bernoulli shift as a factor in each of the three cases:*

(i) $0 < \lambda < \lambda' = 1$,
(ii) $0 < \lambda' < \lambda = 1$,
(iii) $0 < \lambda < \lambda'$.

We recall that i.i.d. Bernoulli shifts are *ergodic* so that the left-shift-invariant sigma-algebra is trivial. It follows from [36, Proposition 4.8] that conservative Bernoulli shifts are also ergodic. Unlike the i.i.d. case, it is not true that the direct product of a conservative Bernoulli shift with itself remains conservative [40, Corollary 6.3]. The *ergodic index* of an ergodic system is the positive integer $k$ such that the $k$-fold direct product of itself remains ergodic, but the $(k+1)$-fold direct product is no longer ergodic. Danilenko and Lemańczyk [7, Question 6] asked what are the possible ergodic indices for type-III$_1$ Bernoulli shifts on two symbols; we will show that every index is possible.

A particular nonsingular Bernoulli shift on two symbols $\{0,1\}$ that we will make use of in answering Danilenko and Lemańczyk's question is given by the half-stationary product measure $\nu^c$ with marginals

$$(1) \qquad \nu_n^c(0) = \frac{1}{2} + \frac{c}{\sqrt{n}} \cdot \mathbf{1}_{[n \geq 1,\, c/\sqrt{n} < 1/2]},$$

where $c > 0$ is a parameter. These shifts were considered by Vaes and Wahl [40, Corollary 6.3], who proved that the shift is of type III$_1$ and ergodic with respect to $\nu^{1/6}$, and for $k \geq 73$ the $k$-fold direct product of the $\nu^{1/6}$ Bernoulli shift is no longer ergodic.

THEOREM 3. *Let $c > 0$. Let $\nu^c$ be the probability measures from (1). There exists $D > 1/6$ such that $(\{0,1\}^{\mathbb{Z}}, \mathcal{B}, \nu^c, T)$ is conservative for all $c < D$ and dissipative for all $c > D$. In addition, if $k \in \mathbb{Z}^+$ and $c \in (D/\sqrt{k+1}, D/\sqrt{k})$, then $(\{0,1\}^{\mathbb{Z}}, \mathcal{B}, \nu^c, T)$ is of ergodic index $k$.*

As a consequence of Theorem 3 and the fact that the ergodic index is nondecreasing under factors, we immediately obtain the following negative result.

COROLLARY 4. *Consider the parameterized measures given in (1). Let $D$ be in as Theorem 3. Let $1 \leq k' < k$ be integers. If $c \in (D/\sqrt{k+1}, D/\sqrt{k})$, but $c' \in (D/\sqrt{k'+1}, D/\sqrt{k'})$, then the Bernoulli shift given by $\nu^{c'}$ is not a factor of $\nu^c$.*

**2. Proof of Theorem 1.** We will sometimes refer to a finite string of symbols as a *block*. If $B = (b_1, \ldots, b_n) = b_1 \cdots b_n$ is a block of binary digits, we let

$$\rho_i(B) := \rho_i(b_1) \cdots \rho_{i+n}(b_n).$$

Two blocks of length 8 that will be important to us are

$$011\ 01\ 011 \quad \text{and} \quad 011\ 10\ 011;$$

here, we inserted spaces to emphasize how the blocks differ.

The following more technical result implies Theorem 1.

THEOREM 5. *A nonsingular Bernoulli shift $\rho = \bigotimes_{i \in \mathbb{Z}} \rho_i$ on two symbols $\{0, 1\}$ satisfying*

$$\sum_{i \in \mathbb{Z}} \left( \frac{\rho_i(01)}{\rho_i(01) + \rho_i(10)} - \frac{1}{2} \right)^2 < \infty, \tag{2}$$

*and*

$$\rho_i(011\ 01\ 011) + \rho_i(011\ 10\ 011) \geq q \tag{3}$$

*for some $q > 0$ for all $i \in \mathbb{Z}$, has an i.i.d. factor.*

*Proof of Theorem 1.* Without loss of generality, we may assume that the Bernoulli shift is on two symbols $\{0,1\}$ instead of a general finite set $A$, since any two-set partition $A = A_0 \cup A_1$ induces a factor map $\phi : A^{\mathbb{Z}} \to \{0,1\}^{\mathbb{Z}}$ given by $\phi(x)_0 = \mathbf{1}_{[x_0 \in A_1]}$, where the push-forward of a product measure on $A^{\mathbb{Z}}$ with the Doeblin condition is again a product measure that has the Doeblin condition on $\{0,1\}^{\mathbb{Z}}$. We will show that nonsingularity and the Doeblin condition imply the two conditions of Theorem 5.

Let $\delta > 0$ be such that for all $i \in \mathbb{Z}$, we have $\delta < \rho_i(0) < 1 - \delta$. Thus for all $i \in \mathbb{Z}$, we have

$$\rho_i(011\ 01\ 011) + \rho_i(011\ 10\ 011) \geq \delta^8.$$

Write $\varepsilon_i := \rho_{i+1}(0) - \rho_i(0)$, so that $\rho_{i+1}(0) = \rho_i(0) + \varepsilon_i$ and

$$\rho_i(1) = 1 - \rho_i(0) = \rho_{i+1}(1) + \varepsilon_i.$$

Elementary manipulations show for $i \in \mathbb{Z}$ that

$$\rho_i(10) = \rho_i(1)\rho_{i+1}(0) = (\rho_{i+1}(1) + \varepsilon_i)(\rho_i(0) + \varepsilon_i)$$
$$= \rho_i(01) + (\rho_i(0) + \rho_{i+1}(1))\varepsilon_i + \varepsilon_i^2 = \rho_i(01) + \varepsilon_i,$$

where the last equality follows from the fact that $\rho_i(0) + \rho_{i+1}(1) = 1 - \varepsilon_i$.

Since the shift is nonsingular, it follows from Kakutani's theorem [18] that

$$\sum_{i \in \mathbb{Z}} \varepsilon_i^2 < \infty;$$

consequently, $\varepsilon_i \to 0$ as $i \to \infty$, and

$$C := \sup_{i \in \mathbb{Z}} \left( \frac{1}{4\rho_i(01) + 2\varepsilon_i} \right)^2 < \infty.$$

Hence,

$$\sum_{i \in \mathbb{Z}} \left( \frac{\rho_i(01)}{\rho_i(01) + \rho_i(10)} - \frac{1}{2} \right)^2 = \sum_{i \in \mathbb{Z}} \left( \frac{\rho_i(01)}{2\rho_i(01) + \varepsilon_i} - \frac{1}{2} \right)^2$$
$$= \sum_{i \in \mathbb{Z}} \left( \frac{1}{4\rho_i(01) + 2\varepsilon_i} \right)^2 \varepsilon_i^2 \leq C \sum_{i \in \mathbb{Z}} \varepsilon_i^2 < \infty. \blacksquare$$

It remains to prove Theorem 5. Let $(\rho_i)_{i \in \mathbb{Z}}$ be a sequence of probability measures on $A$. We say that such a sequence satisfies a *safe zone* condition if there exists a subset $B \subset A$ with at least two elements such that the measures $\rho_i|_B$ are all identical and nonzero, and $\rho_i(B)$ is uniformly bounded away from zero. We proved that any nonsingular Bernoulli shift satisfying a safe zone assumption has an i.i.d. factor [25, Theorems 5 and 7].

Our proof of Theorem 5 builds upon our previous result; we remove the safe zone assumption by using a rudimentary version of Keane and Smorodinsky's [19, 20] marker-filler construction to define a suitable replacement. Previously, Soo and Wilkens [38] used a similar rudimentary marker-filler construction in the i.i.d. setting of Bernoulli actions of a free group to define factors respecting the probabilistic condition of stochastic domination. One of the ingredients in our proof may also be reminiscent of von Neumann's [41] procedure on how to generate a fair coin toss given a possibly bias coin of unknown parameter; see also [33].

Let $X$ be a binary sequence with law $\rho$, a product measure. We say that the integer interval $\{n, n+1, n+2\}$ is a *marker* if

$$X_n X_{n+1} X_{n+2} = 011.$$

Any nonempty interval between two markers is a *filler*; since markers do not overlap, the integers are partitioned into alternating intervals of markers and fillers. We say that a block $B$ *appears* in a block $F$ if $F = KBK'$ for some blocks $K$ and $K'$. Thus the block 011 will not appear in a filler. The *filler measure* on an interval $[k, k']$ is given by conditioning the product measure $\bigotimes_{k \leq i \leq k'} \rho_i$ so that the block 011 does not appear.

We say that a filler $\{n, n+1\}$ is *special* if it is of length 2 and is of the form $(X_n, X_{n+1}) = (1, 0)$ or $(X_n, X_{n+1}) = (0, 1)$. Notice that the filler measure does not require conditioning on a special filler, since the filler is of length 2, and thus markers cannot appear in it.

We say that an integer interval of length 8 is *good* if $X$ restricted to the interval is given by one of the two blocks

011 01 011   or   011 10 011.

Thus, if an interval is good, we know that it contains a special filler.

LEMMA 6. *Let $X$ have law $\rho$ satisfying* (3). *Then there are infinitely many special fillers.*

*Proof.* Partition $\mathbb{Z}$ into intervals of length 8. By (3), each of these intervals is independently good with probability at least $q > 0$. ■

LEMMA 7 (Conditioning). *Let $X$ have law $\rho$ satisfying* (3). *Then the law of $X$ can be sampled by first sampling the markers, and then independently sampling the corresponding filler measures.*

*Proof.* The proof follows from a routine adaptation of [19, Lemma 4], which Keane and Smorodinsky proved for the i.i.d. case. ∎

Given a random sequence $Z \in \{a,b\}^{\mathbb{Z}}$, a *d-equivariant matching scheme* $G$ is an equivariant function of $Z$ such that each integer $m$ with $Z_m = b$ is assigned to an integer $n$ with $Z_n = a$, and each such $n$ has at most $d$ assignments; the assignment function $G$ is equivariant in the sense that if $m$ is assigned to $n$ under $G(Z)$, then $m-1$ is assigned to $n-1$ under $G \circ T(Z)$. Thus a $d$-equivariant matching scheme matches, in a shift-equivariant way, every $b$ to $a$, where each $a$ has at most $d$ partners; see also [25, Section 3].

We will use the same construction used in the proof of [25, Proposition 10], which uses an idea going back to Meshalkin [31]; our description here is adapted from Holroyd and Peres [17]. The *Meshalkin matching scheme* is described inductively as follows. Let $W \in \{a,b\}^{\mathbb{Z}}$ be a random sequence. If $W_{n-1} = b$ and $W_n = a$, then $n-1$ is matched to $n$; that is, in a sequence of $a$'s and $b$'s, we match a $b$ to an $a$ if $b$ is immediately followed by an $a$. In the next step, we disregard all the $b$'s that have been matched, and all the $a$'s that already have been assigned $d$ partners. We repeat this procedure inductively, and we say that it is *successful* if every $b$ is eventually matched to an $a$. In [25, Proposition 10], we proved that $Z$ has a $d$-equivariant matching scheme provided that it is a Bernoulli shift, where the marginal probability of $a$, given by $\mathbb{P}(Z_n = a)$, is uniformly bounded away from $0$. We say that $Z'$ *dominates* $Z$ if $Z_i = a$ implies that $Z'_i = a$.

LEMMA 8 (Monotonicity). *Let $d \in \mathbb{Z}^+$. Let $Z$ and $Z'$ be random processes taking values in $\{a,b\}^{\mathbb{Z}}$. If $Z'$ dominates $Z$ and the Meshalkin $d$-equivariant matching scheme is successful for $Z$, then it is successful for $Z'$.*

*Proof.* Define a random sequence $W$ by setting $W_n = -1$ if $Z'_n = b$ and $W_n = d$ if $Z'_n = a$. Let $m \in \mathbb{Z}$. If $Z'_m = b$, then set
$$R_{Z'} = \inf\{k \geq 1 : W_m + \cdots + W_{m+k} \geq 0\}.$$
Observe that
$$\mathbb{P}(m \text{ is not matched to } m+\ell \text{ for all } \ell \leq k) = \mathbb{P}(R_{Z'} > k).$$
Since $Z'$ dominates $Z$, we have $R_{Z'} \leq R_Z$ and
$$\mathbb{P}(R_{Z'} > k) \leq \mathbb{P}(R_Z > k) \to 0,$$
as desired. ∎

Given an interval $[k, k+n]$, we will sometimes refer to the integer $k$ as the *initial* integer.

PROPOSITION 9 (Matching). *Let $X$ have law $\rho$ satisfying* (3). *There exists an integer $d$ such that the initial integer of each special filler is assigned to at most $d$ other integers in an equivariant way.*

*Proof.* Let $Z' \in \{a,b\}^{\mathbb{Z}}$ be the random sequence where $Z'_n = a$ if $n$ is the initial integer of a special filler of $X$, and $Z_n = b$ otherwise. Note that $Z'$ is not an independent sequence and its corresponding dynamical system, $(\{a,b\}^{\mathbb{Z}}, \mathcal{B}, \mathbb{P}(Z \in \cdot), T)$, is not a Bernoulli system. Consider the indexed partition of the integers of size 8 given by

$$\mathbb{Z} = \bigcup_{n \in \mathbb{Z}} (8n + [0,7]).$$

By (3), each of these intervals is independently good with probability at least $q > 0$. Let $Z \in \{a,b\}^{\mathbb{Z}}$ be the random sequence where $Z_{8n+3} = a$ if the interval $8n + [0,7]$ is good with respect to $X$, and $Z_n = b$ otherwise. Thus $Z_k = a$ only if $k$ an initial integer of special filler of $X$ that occurs in a good interval of this partition, and $Z'$ dominates $Z$. By Lemma 8, it suffices to show that the Meshalkin matching scheme is successful for $Z$. Below is a sample realization of $X$, $Z'$, and $Z$, where $Z$ misses a special filler that is recorded in $Z'$:

| $X$ | $=$ | $\cdots$ | 01101011 | 00000011 | 10011011 | 01110011 | $\cdots$ |
| $Z'$ | $=$ | $\cdots$ | $bbbabbbb$ | $bbbbbbbb$ | $\underline{a}bbbbbbb$ | $bbbabbbb$ | $\cdots$ |
| $Z$ | $=$ | $\cdots$ | $bbbabbbb$ | $bbbbbbbb$ | $bbbbbbbb$ | $bbbabbbb$ | $\cdots$ |

Let $d \geq 8\bigl(1 + \frac{1-q}{q}\bigr)$. Again, define a random sequence $W$ by setting $W_n = -1$ if $Z_n = b$, and $W_n = d$ if $Z_n = a$. Let $m \in \mathbb{Z}$. If $Z_m = b$, then set

$$R_Z = \inf \{k \geq 1 : W_m + \cdots + W_{m+k} \geq 0\}.$$

We will show that $\mathbb{P}(R_Z > k) \to 0$ as $k \to \infty$. Let $E_n = Z_{8n+3}$. Note that $E$ is a Bernoulli shift. Set $Y'_n = -8$ if $E_n = b$, and $Y'_n = d-7$ if $E_n = a$. Then if $m$ lies in the interval $8\ell + [0,7]$, we have

(4) $\qquad W_m + \cdots + W_{m+8k-1} \geq Y'_\ell + Y'_{\ell+1} + \cdots + Y'_{\ell+k-1} - d;$

here we subtract $d$ to account for the possibility that $m$ lies in a good interval, but is to the right of the special filler. Note that the $Y'_n$ are independent, where $\mathbb{P}(Y'_n = d-7) > q$ for all $n \in \mathbb{Z}$. By an elementary version of Strassen's theorem [39], we define an i.i.d. sequence $Y$ with $\mathbb{P}(Y_0 = d-7) = q$ and $\mathbb{P}(Y_0 = -8) = 1-q$ such that for all $n \in \mathbb{Z}$, we have $Y'_n = d-7$ if $Y_n = d-7$. Then

$$\mathbb{E}Y'_n \geq \mathbb{E}Y_n \geq -8(1-q) + \bigl(8\bigl(1+\tfrac{1-q}{q}\bigr) - 7\bigr)q = q > 0.$$

Thus by the law of large numbers it follows that the corresponding partial sums for $Y$ in (4) will become nonnegative almost surely, and thus the same also holds for $Y'$, so that $R_Z$ is finite almost surely. ∎

*Proof of Theorem 5.* Let $X$ have law $\rho$. Let $S = \{n_k, n_{k+1}\}_{k \in \mathbb{Z}}$ be the sequence of special fillers, where we agree that $n_0 \leq 0$ is the largest such integer. Let $Z_k = 1$ if $X_{n_k} X_{n_k+1} = 10$, and $Z_k = 0$ if $X_{n_k} X_{n_k+1} = 01$. By assumption (2) and an application of Katutani's theorem [18, Corollary 1], the law of the sequence $Z$ of the corresponding bits is equivalent to an i.i.d. sequence $Z'$ of fair bits. We remark that the special fillers are conditioned to have the form 10 or 01 and correspondingly the term $\rho_i(01) + \rho_i(10)$ appears in the denominator in (2).

Choose a probability measure $\beta$ on $\{0,1\}$ such that
$$dH(\beta) = -d[\beta(0)\log\beta(0) + \beta(1)\log\beta(1)] = \log 2.$$
Note that $\beta$ will be a biased measure on $\{0,1\}$, and the entropy $H(\beta)$ will be small if $d$ is large. Let $\psi : \{0,1\}^{\mathbb{Z}} \to (\{0,1\}^{d+1})^{\mathbb{Z}}$ be an isomorphism of the uniform product measure to the product measure $(\beta^d)^{\mathbb{Z}}$. Thus, in an equivariant way, we can replace a fair bit by $d$ independent low entropy bits. Note that by Keane and Smorodinsky [20], we may demand the isomorphism $\psi$ is explicit and finitary.

We apply $\psi$ to $Z$ and obtain, in an equivariant manner, a way to associate $d+1$ bits to each special filler. By Proposition 9, we assign, in an equivariant way, bits to all the other integers, retaining one bit for the initial vertex of a special filler, discarding any surplus bits. Thus as a factor of $X$ we obtain a random sequence $W$ of bits.

Notice that the law of $W$ is the push-forward of a function $F(S,Z)$, and the law of $(S,Z)$ is equivalent to that of $(S,Z')$. Since, by construction, $F(S,Z') \stackrel{d}{=} W'$, where $W'$ is an i.i.d. sequence of bits, we conclude that the law of $W$ is equivalent to that of $W'$. ∎

REMARK 10. We recall that the isomorphism $\psi$ in the proof of Theorem 5 is finitary if it is continuous almost surely and has a random coding radius that is finite almost surely. We note that by appealing to Keane and Smorodinsky [20], the factor map given by Theorem 5 is also finitary, since the Meshalkin matching scheme is finitary. See [35] for more information on finitary codes. ◊

### 3. Proof of Theorem 2

**3.1. Essential values.** We give a brief overview of how type-III systems are further parameterized by an additional parameter $\lambda \in [0,1]$; for more details see [34]. Consider a probability space $(\Omega, \mathcal{F}, \mu)$, where $\mu$ is not necessarily a product measure, that is endowed with a group action $G$. We say that $G$ is *nonsingular* if $\mu \circ g$ is equivalent to $\mu$ for all $g \in G$, and *ergodic* if the group-invariant sigma-algebra is trivial. We say that $r \in \mathbb{R}$ is an *essential value* for the group action $G$ if for all $A \in \mathcal{F}$ with positive measure there exist $\varepsilon > 0$ and $g \in G$ such that
$$\mu\left(A \cap g^{-1}(A) \cap \left\{\omega \in \Omega : \left|\log \frac{d(\mu \circ g)}{d\mu}(\omega) - r\right| < \varepsilon\right\}\right) > 0.$$

The *Krieger ratio set*, which consists of the essential values, is a closed subgroup of $\mathbb{R}$. For $\lambda \in (0,1)$, a system is of *type* III$_\lambda$ if the ratio set is $\{n \log \lambda : n \in \mathbb{Z}\}$, and is of *type* III$_1$ if the ratio set is all of $\mathbb{R}$. Thus one way of realizing a type-III$_1$ transformation is to ensure that it contains two rationally independent essential values.

**3.2. Specific type-III Bernoulli shifts.** We define the constructions from [25] that will be used in the proof of Theorem 2. Related constructions are also defined by Berendschot and Vaes [4] and they also constructed type-III$_0$ Bernoulli shifts.

For $n \geq 2$, set
$$a_n := \frac{1}{(n+4)\log(n+4)}.$$

Let $\mathcal{L}(A) = |A|$ denote the Lebesgue measure or length of an interval $A$. Let $\lambda \in (0,1)$. Let $\{A_n\}_{n=2}^\infty$ and $\{B_n\}_{n=2}^\infty$ be decreasing sequences of open intervals of $[0,1]$ satisfying:

(a) For all $n \in \mathbb{N}$, $A_n \cap B_n = \emptyset$.
(b) For all $n \in \mathbb{N}$, $A_{n+1} \subset A_n$ and $B_{n+1} \subset B_n$.
(c) For all $n \in \mathbb{N}$, $|A_n| = a_n = \lambda^{-1}|B_n|$.

Using these sequences we define a sequence of functions $f_n : [0,1] \to \{\lambda^{-1}, 1, \lambda\}$. For all integers $n \leq 1$, set $f_n \equiv 1$. For $n \geq 2$, set

(5) $$f_n(u) := \begin{cases} \lambda, & u \in A_n, \\ 1/\lambda, & u \in B_n, \\ 1, & u \in [0,1] \setminus (A_n \cup B_n). \end{cases}$$

Identify the densities $f_n$ with the associated measures given by
$$E \mapsto \int_E f_n(u)\, du.$$

We proved in [25, Theorem 1] that the Bernoulli shift
$$\left([0,1]^\mathbb{Z}, \mathcal{B}, \bigotimes_{n \in \mathbb{Z}} f_n, T\right)$$
is a nonsingular Bernoulli shift of type III$_\lambda$ that satisfies a *safe zone* condition on a subset of the interval $[0,1]$; namely, there is a set $B \subset [0,1]$ of positive measure such that $\left(\int_B f_n(u)\, du\right)^{-1} f_n \mathbf{1}_B \equiv 1$.

From the definition of the densities in (5), it is clear that the ratio set is a subset of $\{n \log \lambda : n \in \mathbb{Z}\}$, so that one only needs to verify that $\lambda$ is indeed an essential value. This is part of the advantage of working in the continuous setting of densities, as opposed to the case of finitely many symbols, because we can write down candidates for which it is clear that the ratio set is a subset of $\{n \log \lambda : n \in \mathbb{Z}\}$.

REMARK 11. Our proof of [25, Theorem 1] can be summarized as follows:
- Verify nonsingularity and conservativity.
- Verify that $\lambda$ is an essential value of the associated action of the group of all finite permutations on $([0,1]^{\mathbb{Z}}, \mathcal{B}, \bigotimes_{n \in \mathbb{Z}} f_n)$. Note that the action of the group of all finite permutations is ergodic by [25, Lemma 20], as the product measure satisfies the *tameness* condition of Aldous and Pitman [2] for exchangeability.
- Transfer the result on permutations back to the setting of the left-shift. ◊

REMARK 12. It is not surprising that we may *extend* the safe zone by considering the densities

$$g_n = \frac{1}{2}\mathbf{1}_{[-1,0)} + \frac{1}{2}\mathbf{1}_{[0,1]} f_n,$$

where $f_n$ is given in (5); our proof of [25, Theorem 1] implies that the Bernoulli shift $([0,1]^{\mathbb{Z}}, \mathcal{B}, \bigotimes_{n \in \mathbb{Z}} g_n, T)$ is still of type III$_\lambda$. ◊

REMARK 13. Consider the type-III$_1$ Bernoulli shift given in the following way. Let $L, \lambda \in (0,1)$. Let $\rho$ and $\nu$ correspond to the type-III$_\lambda$ and III$_L$ Bernoulli shifts as given above, where the marginals are defined on disjoint sets $[-1, 0)$ and $[0, 1]$. If $\log \lambda$ and $\log L$ are rationally independent, it is easy to see that any nontrivial convex combination of $\rho$ and $\nu$ corresponds to a type-III$_1$ Bernoulli shift. See also [25, Example 33]. ◊

Our proof of Theorem 2 consists in considering the type-III$_1$ system from Remark 13, and *erasing* one the original type-III shifts by replacing it with the uniform distribution.

*Proof of Theorem 2(ii).* Consider the type-III$_1$ system from Remark 13, with $\log \lambda$ and $\log L$ rationally independent and the measure

$$\mu = \bigotimes_{i \in \mathbb{Z}} \left(\tfrac{1}{2}\rho_i + \tfrac{1}{2}\nu_i\right).$$

Thus if the random sequence $Z$ has law $\mu$, then $Z$ can be sampled by independently sampling $Z_i$ from $\rho_i$ or $\nu_i$ with equal probability. In addition, note that $\rho_i$ and $\nu_i$ have disjoint supports. We say that $i$ is a $\nu$-*index* if $Z_i < 0$. Let $(n_k)_{k \in \mathbb{Z}}$ be an enumeration of the $\nu$-indices where $n_k < n_{k+1}$ and we agree that $n_0 \leq 0$ is the largest such integer.

Let $\nu$ satisfy the safe zone condition with a set $B \subset [-1, 0)$. We say a $\nu$-index $i$ is *special* if $Z_i \in B$. Thus conditional on $i$ being a special index, $Z_i$ is uniformly distributed on $B$. Hence at each special index $i$, as a deterministic function of $Z_i$, we can produce an infinite sequence of random variables that are uniformly distributed on $[-1, 0)$. By a routine argument [25, Theorem 5], we distribute these uniform random variables in an equivariant way to each $\nu$-index. Thus as a factor of $Z$, we can produce an independent sequence $Z'$, where $Z'_i$ is uniformly distributed on $[-1, 0)$ if $i$ is a $\nu$-index, and $Z'_i = Z_i$ otherwise.

Hence
$$\mu' = \bigotimes_{i \in \mathbb{Z}} \left(\tfrac{1}{2}\mathcal{L}|_{[-1,0)} + \tfrac{1}{2}\rho_i\right)$$
is a factor of $\mu$ and it follows from Remark 12 that the corresponding Bernoulli shift is of type III$_\lambda$. ∎

**3.3. Piecewise linear transformations.** Our proof of part (iii) of Theorem 2 will be obtained by applying a single piecewise linear transformation $h$ to each coordinate so that our factor map will be of the form
$$[\phi(x)]_i = h(x_i).$$
We will make use of the following special case of an elementary change of variables formula [14, equation (16), p. 112].

LEMMA 14. *Let $U$ be a real-valued random variable with a probability density function $f_U$ and $h$ be piecewise linear and finite-to-one. Then the density of the random variable $V = h(U)$ is given by*
$$f_V(v) = \sum_{u \in h^{-1}(v)} f_U(u)(|h'(v)|)^{-1},$$
*where we set $f_V(h(u)) = 0$ if $h$ is not differentiable at $u$.*

With Lemma 14, we will be able to modify the densities from Section 3.2 so that the resulting densities have the essential value(s) we desire.

*Proof of Theorem 2(iii).* Let $0 < \lambda < \lambda' < 1$. Set
$$(6) \qquad p := \frac{\lambda' - \lambda}{1 - \lambda'}.$$
Consider the densities $f_n$ in Section 3.2 as expressed in (5), so that the corresponding Bernoulli shift is of type III$_\lambda$. With minor modifications, such as re-indexing the sequence $a_n$, we assume that
$$a_1 + pa_1 < 1/2.$$
We will also write $a_n = 0$ for all $n \leq 0$. Define $h : [0,1] \to [0,1]$ by
$$h(x) = \begin{cases} (x - a_1)/p, & a_1 < x < a_1 + pa_1, \\ (1 - \lambda a_1) + \lambda \frac{(1-\lambda a_1)-x}{p}, & 1 - \lambda a_1 - pa_1 < x < 1 - \lambda a_1, \\ x, & \text{otherwise.} \end{cases}$$

Set $\phi : [0,1]^{\mathbb{Z}} \to [0,1]^{\mathbb{Z}}$ via
$$[\phi(x)]_n = h(x_n).$$
By Lemma 14, and routine calculations, the push-forward of the type-III$_\lambda$ product measures satisfies
$$\left( \bigotimes_{n \in \mathbb{Z}} f_n \right) \circ \phi^{-1} =: \bigotimes_{n \in \mathbb{Z}} g_n$$
with
$$g_n(x) = \begin{cases} \lambda + p, & 0 < x < a_n, \\ 1 + p, & a_n < x < a_1, \\ 1, & a_1 + pa_1 < x < 1 - (\lambda + p)a_1, \\ 1 + p/\lambda, & 1 - \lambda a_1 \leq x < 1 - \lambda a_n, \\ (1+p)/\lambda, & 1 - \lambda a_n < x < 1, \\ 0, & \text{otherwise}, \end{cases}$$
for all $n \geq 1$; when $n \leq 0$, we have
$$g_n(x) = \begin{cases} 1 + p, & 0 < x < a_1, \\ 1, & a_1 + pa_1 < x < 1 - (\lambda + p)a_1, \\ 1 + p/\lambda, & 1 - \lambda a_1 \leq x < 1, \\ 0, & \text{otherwise}. \end{cases}$$
Note that the closure of the support of $g_n$ does not depend on $n$. We write $I \subset [0,1]$ for the closure of the support of the $g_n$'s.

We will now argue that the Bernoulli shift corresponding to the product measure $\bigotimes_{n \in \mathbb{Z}} g_n$ is of type III$_{\lambda'}$. Again, this newly constructed Bernoulli shift is a factor of one that is conservative. We will show that the only possible essential value is $\log \lambda'$, from which the same reasoning as in our proof of Theorem 2(i) implies that the new Bernoulli shift is of type III$_{\lambda'}$.

From our expressions for $g_n$, we deduce for every $n \in \mathbb{Z}$ and $v \in I$ that
$$\frac{g_{n-1}(v)}{g_n(v)} \in \left\{ \frac{\lambda + p}{1 + p}, 1, \frac{\lambda^{-1}(1+p)}{1 + p/\lambda} \right\}.$$
Our initial choice of $p$ in (6) was such that
$$\frac{\lambda + p}{1 + p} = \lambda' \quad \text{and} \quad \frac{\lambda^{-1}(1+p)}{1 + p/\lambda} = \frac{1 + p}{\lambda + p} = \frac{1}{\lambda'}.$$
Hence the Krieger ratio set of the new Bernoulli shift is contained in $(\log \lambda')\mathbb{Z}$. A routine variation of the argument given in [25, Theorem 23] tells us that $\log \lambda'$ is an essential value with respect to the action of the group of finite permutations. Finally, by [3, Section 4] we find that the result regarding permutations can be exchanged for the desired result with respect to the left-shift. ∎

*Proof of Theorem 2(i).* Consider the type-III$_\lambda$ Bernoulli shift defined as follows. Let $\lambda \in (0,1)$. Let $\rho$ and $\nu$ correspond to the type-III$_\lambda$ Bernoulli shifts as given in (5), where the marginals are defined on the disjoint sets $[-1, 0)$ and $[0, 1]$. Then the measure

$$\mu = \bigotimes_{i \in \mathbb{Z}} \left(\tfrac{1}{2}\rho_i + \tfrac{1}{2}\nu_i\right)$$

is also of type III$_\lambda$.

Choose $\lambda' > \lambda$ with $\log \lambda'$ and $\log \lambda$ rationally independent. Let $Z = (Z_i)_{i \in \mathbb{Z}}$ have law $\mu$. As in the proof of Theorem 2(ii), we will only modify $Z$ when $Z_i < 0$. We obtain a random sequence $Z'$, by setting $Z'_i = Z_i$ if $Z_i \geq 0$, and by applying the coordinatewise factor map from Theorem 2(iii), which takes $\lambda$ to $\lambda'$, to $Z_i$ if $Z_i < 0$. Let $\mu'$ be the law of $Z'$. Clearly, $\mu'$ is a Bernoulli shift and is obtained as a factor from $\mu$. By Remark 13, the measure $\mu'$ is of type III$_1$. ∎

**4. Some tools for the proof of Theorem 3.** We saw in the proof of Theorem 2 that we would deduce that a system was conservative simply because it was a factor of a conservative system. Let $(\Omega, \mathcal{F}, \mu, S)$ and $(\Omega', \mathcal{C}, \eta, R)$ be nonsingular systems. We will often refer to these systems by the transformations $S$ and $R$ or their measures $\mu$ and $\eta$. We say that $\eta$ is a *factor* of $\mu$ if there exists a measurable function, a *factor map*, $\phi : \Omega \to \Omega'$ such that $\phi$ is equivariant so that $\phi \circ S = R \circ \phi$, and the push-forward of $\mu$ under $\phi$ is equivalent to $\eta$; moreover, we say that $\eta$ and $\mu$ are *isomorphic* if $\phi^{-1}$ serves as a factor map from $\mu$ to $\eta$.

Note that in the case of an i.i.d. Bernoulli shift, the $k$-fold direct product of the left-shift is isomorphic to the $k$-fold *composition* of the left-shift. A version of this elementary fact (Theorem 20) for nonsingular Bernoulli shifts will be used in our proof of Theorem 3.

Also consider the case where $(\Omega, \mathcal{F}, S) = (\Omega', \mathcal{C}, R)$. If $\mu$ and $\eta$ are equivalent measures on $\Omega$, then we will write $\mu \sim \nu$; in this case, the identity map is an isomorphism between $(\Omega, \mathcal{F}, \mu, S)$ and $(\Omega, \mathcal{F}, \eta, S)$.

Given two systems it is often easier to study them if they are factors of some constructed larger system that is not merely the direct product of the two; such techniques are loosely associated with *coupling* in probability theory [9, 30] and also *joining* in ergodic theory [12, 13].

**4.1. Some operations on product measures.** Given a product measure $\mu$ on $A^{\mathbb{Z}}$ and $0 < p < 1$ there are several natural operations on $\mu$ which are done by tossing a coin infinitely many times and then using the outcome as an indication of whether to use the data of $\mu$ or some external source.

Let $\mu = \bigotimes_{i \in \mathbb{Z}} \mu_i$ be a product measure on $A^{\mathbb{Z}}$. Let $p \in (0,1)$ and $\{H, T\}$ be a set of two elements. We let the probability vector $(p, 1-p)$ also denote the probability measure on $\{H, T\}$ giving mass $p$ to H and $1 - p$ to T. Similarly, $(p, 1 - p)$ is a probability measure on the set $\{0, 1\}$, giving mass $p$ to the element $0$. Let $\alpha$ be a probability measure on $A$. Consider the product space $A^{\mathbb{Z}} \times \{H, T\}^{\mathbb{Z}} \times A^{\mathbb{Z}}$ endowed with the product measure $\mu \otimes (p, 1 - p)^{\mathbb{Z}} \otimes \alpha^{\mathbb{Z}}$ and the shift $T \times T \times T$. Consider the equivariant maps $\Theta : A^{\mathbb{Z}} \times \{H, T\}^{\mathbb{Z}} \times A^{\mathbb{Z}} \to A^{\mathbb{Z}} \times \{H, T\}^{\mathbb{Z}}$ and $\Phi : A^{\mathbb{Z}} \times \{H, T\}^{\mathbb{Z}} \times A^{\mathbb{Z}} \to A^{\mathbb{Z}}$ given by

$$\Theta(x, y, z)_j = \begin{cases} (x_j, H), & y_j = H, \\ (z_j, T), & y_j = T, \end{cases} \quad \Phi(x, y, z)_j = \begin{cases} x_j, & y_j = H, \\ z_j, & y_j = T. \end{cases}$$

Thus, with probability $p$ we stay with the measure $\mu$, and with probability $1 - p$ we choose an output from the measure $\alpha$. The *random insertion operation on $\mu$ of parameters $p$ and $\alpha$* is given by the push-forward

$$\mathrm{RI}(\mu, p, \alpha) := (\mu \otimes (p, 1 - p)^{\mathbb{Z}} \otimes \alpha^{\mathbb{Z}}) \circ \Theta^{-1}.$$

The *randomized product measure of parameters $p$ and $\alpha$* is given by the push-forward

$$\mathrm{RPM}(\mu, p, \alpha) := (\mu \otimes (p, 1 - p)^{\mathbb{Z}} \otimes \alpha^{\mathbb{Z}}) \circ \Phi^{-1}.$$

Both these measures are factors of the original triple-product measure, and the randomized product measure is obtained as a measure-preserving factor of the random insertion operation. Thus, if $\mu$ is nonsingular with respect to the shift on $A^{\mathbb{Z}}$, then both $\mathrm{RI}(\mu, p, \alpha)$ and $\mathrm{RPM}(\mu, p, \alpha)$ are nonsingular with respect to the shift on their corresponding product spaces.

EXAMPLE 15. Let $p \in (0, 1)$ and $a_n \to 0$ as $n \to \infty$. Consider the product measure $\mu = \bigotimes_{n \in \mathbb{Z}} \mu_n$ on $\{0, 1\}^{\mathbb{Z}}$ with marginals $\mu_n(0) = p$ if $p + a_n < 0$ or $p + a_n > 1$, and otherwise

$$\mu_n(0) = 1 - \mu_n(1) = p + a_n.$$

Let $0 < q < 1$. If $\hat{\mu} = \mathrm{RPM}(\mu, q, (p, 1 - p))$, then $\hat{\mu}$ is a product measure on $\{0, 1\}^{\mathbb{Z}}$ with $\hat{\mu}_n(0) = p$ if $p + a_n < 0$ or $p + a_n > 1$, and otherwise

$$\hat{\mu}_n(0) = 1 - \hat{\mu}_n(1) = p + q a_n. \diamond$$

Example 15 tells us that the two measures of the form (1) may be coupled together in the same probability space, where one is a factor of the other with additional randomization.

**4.2. Phase transition.** It was proved, in increasing levels of generality, that a nonsingular Bernoulli shift is conservative or totally dissipative [5, 24]; furthermore, in the conservative case it must be *weakly mixing* so that its direct product with every ergodic probability-preserving system remains ergodic [5, 6].

We will prove the existence of the following phase transition from dissipativity to conservativity on the type-III$_1$ Bernoulli shift given by (1); this phase transition will be important in our proof of Theorem 3.

Let $p \in (0,1)$. Consider the following class of bi-infinite sequences. Let $a_n \to 0$ as $|n| \to \infty$. Let $\mu = \bigotimes_{n \in \mathbb{Z}} \mu_n$ be the product measure on $\{0,1\}^{\mathbb{Z}}$ with marginals satisfying

$$(7) \qquad \mu_n(0) = 1 - \mu_n(1) = p + a_n,$$

except on possibly finitely many coordinates, where we set $\mu_n(0) = p$ if $p + a_n < 0$ or $p + a_n > 1$. In this case, there exists $0 < q < 1$ such that for all $n \in \mathbb{Z}$, we have $q < \mu_n(0) < 1-q$ or in other words $\mu$ satisfies the Doeblin condition. An application of Kakutani's theorem [18] shows that the Bernoulli shift is $\mu$-nonsingular if and only if

$$(8) \qquad \sum_{n \in \mathbb{Z}} (a_n - a_{n-1})^2 < \infty.$$

If $(a_n)_{n \in \mathbb{Z}}$ satisfies the above properties, then we say it is a *nonsingular sequence*.

Let $p \in (0,1)$. For a fixed nonsingular sequence $(a_n)_{n \in \mathbb{Z}}$ consider the linear transformation $a_n \mapsto c a_n$ for $c \in (0, \infty)$. If $\mu$ is the associated product measure on $\{0,1\}^{\mathbb{Z}}$ given by (7) we will write $\mu^{(p,c)}$ for the product measure on $\{0,1\}^{\mathbb{Z}}$ with the marginals $\mu_n^{(p,c)}(0) = p + c a_n$ where again we define the marginals to be $\mu_n^{(p,c)}(0) = p$ at the finitely many integers where $p + c a_n < 0$ or $p + c a_n > 1$.

THEOREM 16. *Let $p \in (0, 1/2]$. Let $(a_n)_{n \in \mathbb{Z}}$ be a nonsingular sequence. There exists $c_0(p) \in [0, \infty]$ such that the nonsingular Bernoulli shift given by $(\{0,1\}^{\mathbb{Z}}, \mathcal{B}, \mu^{(p,c)}, T)$ is totally dissipative for every $c > c_0(p)$ and conservative and ergodic for every $c < c_0(p)$. Furthermore, for all $0 < p < q \leq 1/2$, we have the lower bound $c_0(p) \geq (p/q) c_0(q)$.*

*Proof.* As a Bernoulli shift is either totally dissipative or conservative and in the latter case it is ergodic, in order to show the existence of $c_0(p)$ it suffices to show that if $d < c$ and (the Bernoulli shift corresponding to the product measure) $\mu^{(p,c)}$ is conservative, then $\mu^{(p,d)}$ is conservative. Define $c_0(p)$ as the supremum of all $c$ for which $\mu^{(p,c)}$ is conservative.

Let $d < c$ be such that $(\{0,1\}^{\mathbb{Z}}, \mathcal{B}, \mu^{(p,c)}, T)$ is conservative; this system is ergodic and weakly mixing [6]. Since i.i.d. Bernoulli shifts are probability-preserving and ergodic, the measure

$$\mu^{(p,c)} \otimes (d/c, 1 - d/c)^{\mathbb{Z}} \otimes (p, 1-p)^{\mathbb{Z}} \quad \text{on } \{0,1\}^{\mathbb{Z}} \times \{\text{H}, \text{T}\}^{\mathbb{Z}} \times \{0,1\}^{\mathbb{Z}}$$

is ergodic with respect to the product shift $T \times T \times T$. By Example 15, the marginals of $\text{RPM}(\mu^{(p,c)}, d/c, (p, 1-p))$ and $\mu^{(p,c)}$ are the same up to finitely many $n \in \mathbb{N}$, namely those positive integers $n$ where $p + ca_n \notin (0,1)$, but $p + da_n \in (0,1)$. Thus
$$\mu^{(p,d)} \sim \mathrm{RPM}(\mu^{(p,c)}, d/c, (p, 1-p)),$$
and moreover $\mu^{(p,d)}$ is isomorphic to a factor of an ergodic system, hence remains ergodic.

Finally, we note that an ergodic invertible system is conservative and therefore the shift with respect to $\mu^{(p,d)}$ is conservative, finishing the proof of existence of $c_0(p)$.

For the final remark, let $p < q \leq 1/2$ and observe that for the degenerate distribution $(0,1)$ on $\{0,1\}$, it is easy to see that
$$\mu^{(p,pc/q)} \sim \mathrm{RPM}(\mu^{(q,c)}, p/q, (0,1)).$$
Thus if $c < c_0(q)$, then the measure $\mu^{(p,pc/q)}$ is isomorphic to a factor of the conservative system
$$(\{0,1\}^{\mathbb{Z}} \times \{0,1\}^{\mathbb{Z}}, \mathcal{B} \otimes \mathcal{B}, \mu^{(q,c)} \times (p/q, 1-p/q)^{\mathbb{Z}}, T \times T),$$
and is conservative, so that $c_0(p) \geq (p/q)c_0(q)$. ∎

REMARK 17. We have not investigated the ergodicity of the system *at* the critical value $c_0(p)$, even for specific nonsingular sequences. ◊

REMARK 18. Previously, we argued that the *upper* bound $c_0(p) \leq \frac{p}{q}c_0(q)$ also holds, but it is possible to show that this is not true using Theorem 3. We thank the referee for pointing out a flaw in our previous argument. These inequalities are no longer required for the proof of Theorem 3. ◊

**4.3. Speedups.** Let $T$ be the left-shift. Given $k \in \mathbb{Z}^+$, the $k$-fold composition of $T$ with itself given by $T^k = T^{k-1} \circ T$ will be referred to as the *k-speedup* of $T$. The following proposition which connects the conservativity of $T$ with its speedups follows from a result of Halmos [15].

PROPOSITION 19. *Let $A$ be a finite set and $(A^{\mathbb{Z}}, \mathcal{B}, \mu, T)$ be a nonsingular Bernoulli shift. If $\mu$ satisfies the Doeblin condition, then $\mu$ is conservative with respect to $T$ if and only if for all $k \in \mathbb{Z}^+$ the speedup $(A^{\mathbb{Z}}, \mathcal{B}, \mu, T^k)$ is ergodic.*

*Proof.* Since $\mu$ satisfies the Doeblin condition and $A$ is finite, either $\mu$ is conservative or it is totally dissipative [24]. Fix $k \in \mathbb{Z}^+$. It follows from a result of Halmos [1, Corollary 1.1.4] that $\mu$ is conservative with respect to $T$ if and only if $\mu$ is conservative with respect to the speedup $T^k$. It remains to argue that the speedup is ergodic.

Let $\eta$ be the product measure on $(\{0,1\}^k)^{\mathbb{Z}}$ with marginals
$$(9) \qquad \eta_n = \bigotimes_{i=0}^{k-1} \mu_{kn+i}.$$

The map $\zeta : \{0,1\}^{\mathbb{Z}} \to \left(\{0,1\}^k\right)^{\mathbb{Z}}$ given by
$$\zeta(x)_n = (x_{kn}, \ldots, x_{kn+(k-1)}) \tag{10}$$
is an isomorphism of $T^k$ on $\{0,1\}^{\mathbb{Z}}$ and the left-shift on $(\{0,1\}^k)^{\mathbb{Z}}$ and $\eta = \mu \circ \zeta^{-1}$. Thus the speedup is isomorphic to a Bernoulli shift and conservativity is equivalent to ergodicity. Therefore the speedup $(A^{\mathbb{Z}}, \mathcal{B}, \mu, T^k)$ is ergodic if and only if the Bernoulli shift $(A^{\mathbb{Z}}, \mathcal{B}, \mu, T)$ is conservative. ∎

Proposition 19 together with the following isomorphism will allow us to analyze the ergodicity of the $k$-fold direct product.

THEOREM 20. *Let $p \in (0, 1/2]$. Let $(a_n)_{n \in \mathbb{Z}}$ be a nonsingular sequence and $\mu$ be the product measure on $\{0,1\}^{\mathbb{Z}}$ with marginals*
$$\mu_n(0) := p + a_n.$$
*Then for all $k \in \mathbb{Z}^+$, the speedup $(\{0,1\}^{\mathbb{Z}}, \mathcal{B}, \mu, T^k)$ and the direct product $((\{0,1\}^{\mathbb{Z}})^k, \mathcal{B}^{\otimes k}, \gamma[k]^{\otimes k}, T^{\otimes k})$ are isomorphic, where $\gamma[k]$ is the product measure with marginals given by*
$$\gamma[k]_n(0) = p + a_{kn}.$$

*Proof.* Fix $k \in \mathbb{Z}^+$. We already know from the proof of Proposition 19 that $\zeta$ given by (10) is an isomorphism of the speedup $T^k$ on $\{0,1\}^{\mathbb{Z}}$ and the left-shift on $(\{0,1\}^k)^{\mathbb{Z}}$ endowed with the product measure $\eta$ given by (9).

The map $\pi : (\{0,1\}^{\mathbb{Z}})^k \to (\{0,1\}^k)^{\mathbb{Z}}$ given by
$$\pi(x^1, \ldots, x^k)_n = (x^1_n, \ldots, x^k_n)$$
is an isomorphism between $\gamma[k]$ and the product measure $\kappa$ on $(\{0,1\}^k)^{\mathbb{Z}}$ with marginals given by
$$\kappa_n(b_1, \ldots, b_k) = \bigotimes_{i=1}^{k} \gamma[k]_n(b_i),$$
for $(b_1, \ldots, b_k) \in \{0,1\}^k$. It remains to show that $\kappa$ and $\eta$ are isomorphic, which we will accomplish by the identity map and Kakutani's theorem [18]. The following is a diagram depicting the isomorphisms:

$$\begin{array}{ccc}
((\{0,1\}^k)^{\mathbb{Z}}, \kappa, T) & \xleftarrow{\pi} & ((\{0,1\}^{\mathbb{Z}})^k, \gamma[k]^{\otimes k}, T^{\otimes k}) \\
{\scriptstyle \mathrm{id}} \downarrow & & \\
((\{0,1\}^k)^{\mathbb{Z}}, \eta, T) & \xleftarrow{\zeta} & (\{0,1\}^{\mathbb{Z}}, \mu, T^k)
\end{array}$$

Since $\kappa$ and $\eta$ satisfy the Doeblin condition, it follows from Kakutani's theorem [18] that the measures $\kappa$ and $\eta$ are equivalent if and only if
$$\sum_{n \in \mathbb{Z}} \sum_{B \in \{0,1\}^k} (\eta_n(B) - \kappa_n(B))^2 < \infty, \tag{11}$$

or equivalently, for all $B \in \{0,1\}^k$,
$$\sum_{n \in \mathbb{Z}} (\eta_n(B) - \kappa_n(B))^2 < \infty.$$

Fix $B = (b_1, \ldots, b_k) \in \{0,1\}^k$. Set $\alpha(n) := |\eta_n(B) - \kappa_n(B)|^2$. For all $n \in \mathbb{Z}$, we have

$$\begin{aligned}
(12) \quad \alpha(n) &= \Big| \prod_{j=1}^{k} \mu_{kn+j-1}(b_j) - \prod_{j=1}^{k} \mu_{kn}(b_j) \Big|^2 \\
&= \Big| \sum_{\ell=1}^{k} \Big[ (\mu_{kn+\ell-1}(b_\ell) - \mu_{kn}(b_\ell)) \prod_{j=1}^{\ell-1} \mu_{kn+j-1}(b_j) \prod_{j=\ell+1}^{k} \mu_{kn}(b_j) \Big] \Big|^2 \\
&\leq \Big( \sum_{\ell=1}^{k} |\mu_{kn+\ell-1}(b_\ell) - \mu_{kn}(b_\ell)| \Big)^2 = \Big( \sum_{\ell=1}^{k} |a_{kn+\ell-1} - a_{kn}| \Big)^2 \\
&\leq k^2 \sum_{\ell=1}^{k} |a_{kn+\ell-1} - a_{kn}|^2.
\end{aligned}$$

Since $(a_n)_{n \in \mathbb{Z}}$ is a nonsingular sequence, for all $1 \leq \ell \leq k$, the speedup $T^\ell$ is nonsingular with respect to $\mu$, and again by Kakutani's theorem [18], we have
$$\sum_{n \in \mathbb{Z}} (a_n - a_{n-\ell})^2 < \infty.$$

The right hand side of the inequality (12) is summable over the indices $n$ and (11) holds, yielding the desired equivalence $\eta \sim \kappa$. ∎

**5. The proof of Theorem 3.** We have now assembled all the ingredients necessary for the proof of Theorem 3. Our proof consists in relating the $k$-fold direct product to the $k$-speedup, from which a scaling argument can be applied to the family of measures given by (1).

*Proof of Theorem 3.* As in (1), let $\nu^c$ be the product measure on $\{0,1\}^{\mathbb{Z}}$ with marginals
$$\nu_n^c(0) := \frac{1}{2} + \frac{c}{\sqrt{n}} \cdot \mathbf{1}_{[n \geq 1,\, c/\sqrt{n} < 1/2]}.$$

We will start by showing the existence of $D \in (0, \infty)$ which is a phase transition between conservative and dissipative behaviour. Vaes and Wahl [40, Corollary 6.3] proved that for all $c \leq 1/6$, the measure $\nu^c$ is conservative. Let $a_n(c) := \frac{c}{\sqrt{n}} \cdot \mathbf{1}_{[n \geq 1,\, c/\sqrt{n} < 1/2]}$. In addition, Vaes and Wahl show that

$$\int \sqrt{\frac{d\nu^c \circ T^k}{d\nu^c}} \, d\nu^c \leq \exp\Big( -\frac{1}{2} \sum_{n \in \mathbb{Z}} (a_{n-k}(c) - a_n(c))^2 \Big).$$

For all $k \in \mathbb{Z}^+$, we have

$$\sum_{n \in \mathbb{Z}} (a_{n-k}(c) - a_n(c))^2 \geq c^2 \sum_{\ell=1}^{k} \frac{1}{\ell} \geq c^2(\log k + 1) \quad \text{as } k \to \infty,$$

since 1 is greater than the Euler–Mascheroni constant [29]. Thus if $c > \sqrt{2}$, then

$$\sum_{k=1}^{\infty} \int \sqrt{\frac{d\nu^c \circ T^k}{d\nu^c}} \, d\nu^c < \infty.$$

It follows from [23] that $\mu^c$ is dissipative for all $c > \sqrt{2}$.

Hence by Theorem 16 there exists $1/6 \leq D \leq \sqrt{2}$ such that for all $c < D$, the measure $\nu^c$ is conservative, and for all $c > D$, the measure $\nu^c$ is dissipative.

Now we will apply Theorem 20 to show that this family of measures exhibits every ergodic index. Let $k \in \mathbb{Z}^+$. Since for all large $n \in \mathbb{Z}^+$, we have $a_{kn}(c) = a_n(c)/\sqrt{k}$, we deduce that the product measure $\gamma^c[k]$ with marginals

$$(\gamma^c[k])_n(0) := 1/2 + a_{kn}(c),$$

is equivalent to $\nu^{c/\sqrt{k}}$. By Theorem 20, for all $k \in \mathbb{N}$ and $c > 0$, we have the isomorphism

(13) $\qquad (({\{0,1\}}^{\mathbb{Z}})^k, \mathcal{B}^{\otimes k}, (\nu^c)^{\otimes k}, T^{\otimes k}) \cong (\{0,1\}^{\mathbb{Z}}, \mathcal{B}, \nu^{c\sqrt{k}}, T^k).$

Let $c \in (D/\sqrt{k+1}, D/\sqrt{k})$. Since $c\sqrt{k} < D$, the measure $\nu^{c\sqrt{k}}$ is conservative. By Proposition 19, the $k$-speedup $(\{0,1\}^{\mathbb{Z}}, \mathcal{B}, \nu^{c\sqrt{k}}, T^k)$ is ergodic. Thus the isomorphism (13) implies that the $k$-fold direct product given by $(({\{0,1\}}^{\mathbb{Z}})^k, \mathcal{B}^{\otimes k}, (\nu^c)^{\otimes k}, T^{\otimes k})$ is ergodic.

On the other hand, since $c\sqrt{k+1} > D$, the measure $\nu^{c\sqrt{k+1}}$ is dissipative. Similarly, again by Proposition 19, the $(k+1)$-speedup is dissipative and the isomorphism (13) implies that $(k+1)$-fold direct product is dissipative and hence not ergodic. We conclude that the ergodic index of the shift with respect to $\nu^c$ is precisely $k$. ∎

**Acknowledgments.** We thank Mike Hochman for asking whether any known type-III$_1$ binary Bernoulli shift has an i.i.d. factor. We thank the referee for very helpful comments and questions.

Zemer Kosloff is funded in part by ISF grant No. 1570/17.

Zemer Kosloff  
Einstein Institute of Mathematics  
Hebrew University of Jerusalem  
Edmund J. Safra Campus  
Givat Ram, Jerusalem, 9190401, Israel  
E-mail: zemer.kosloff@mail.huji.ac.il  
http://math.huji.ac.il/~zemkos/

Terry Soo  
Department of Statistical Science  
University College London  
Gower Street  
London WC1E 6BT, United Kingdom  
E-mail: math@terrysoo.com  
http://www.terrysoo.com